\documentclass[11pt, a4paper]{article}
\usepackage{amsfonts,enumerate,amsmath,amsfonts,xspace,array,delarray,theorem,amssymb,epsfig}
\usepackage{graphicx,color}

\theoremstyle{break}\theorembodyfont{\it}\theoremheaderfont{\normalfont\bfseries
}

\newtheorem{theo}{Theorem}

\newtheorem{lem}[theo]{Lemma}

\newtheorem{prop}[theo]{Proposition}

\newtheorem{coro}[theo]{Corollary}

\newenvironment{proof}{\noindent{\bf Proof: }}
                {\leavevmode\unskip\nobreak\hskip2em plus1fill
                $\scriptstyle\bullet$\vskip\theorempostskipamount\par}

\def\LP{Littlewood--Paley\xspace}
\def\dis{\displaystyle}
\def\gh{Heisenberg group\xspace}

\let\lt=<
\let\gt=>

\def\R{{\mathbb R}}
\def\C{{\mathbb C}}
\def\N{{\mathbb N}}
\def\Z{{\mathbb Z}}

\def\H{{\mathbb H}}
\def\LL{{\cal L}}
\def\GG{{\cal G}}

\def\hn{\H_n}

\def\met#1{{\frac{#1}{2}}}
\def\lmnu{L_m^{(n-1)}}

\def\lp#1{L^{#1}(\H_n)}

\def\cb#1#2{\left(\begin{array}{c}#1\\#2\end{array}\right)}
\def\coefbi{\cb{m+n-1}{m}}
\def\ppp#1#2{#1_1,\ldots,#1_{#2}}
\def\ddd#1{\frac{\partial}{\partial #1}}
\def\sjn{\sum_{j=1}^n}
\def\sh{{\cal S}(\hn)}

\def\ir#1{\int_{\R}#1|\lambda|^n\,d\lambda}
\def\cn{\frac{2^{n-1}}{\pi^{n+1}}}
\def\lrad#1{L^{#1}_{\mathit rad}(\hn)}
\def\srad{{\cal S}_{\mathit rad}(\hn)}
\def\okl{\omega_{m,\lambda}}
\def\ih{\int_{\H_n}}
\def\idd#1{\cn\sum_{m=0}^{+\infty}\cb{m+n-1}{m}\ir{#1}}
\def\smo{\sum_{m=0}^{+\infty}}
\def\kl{(m,\lambda)}
\def\mez{\met{1}}

\def\sih{{\cal S}'(\hn)}

\def\be{\begin{equation}}
\def\beq{\begin{eqnarray*}}
\def\eeq{\end{eqnarray*}}

\def\spc#1#2#3#4#5{{L}_{#2}^{#1}(\dot B ^{{#3},{#4}}_{#5}(\LL))}

\def\sp#1#2#3{{L}_{#2}^{#1}\left(#3\right )}
\def\bpd#1#2#3{\dot B ^{{#1},{#2}}_{#3}(\Delta)}
\def\bpl#1#2#3{\dot B ^{{#1},{#2}}_{#3}(\LL)}
\def\bp#1#2#3{\dot B ^{{#1},{#2}}_{#3}(L)}

\def\fj{\varphi_j}
\def\dj{2^{-2j}}
\def\sjz{\sum_{j\in\Z}}

\def\zs{(z,s)}
\def\ddj{2^{2j}}
\def\dnj{2^{Nj}}

\def\t{t^{-\mez}}
\def\tm{|t|^{-\mez}}
\def\pair#1#2{\langle#1,#2\rangle}
\def\RL#1{R(2^{-2#1}(4n\lambda+\lambda^2))}
\def\bL#1#2#3{\dot B ^{{#1},{#2}}_{#3}(\LL)}
\def\bD#1#2#3{\dot B ^{{#1},{#2}}_{#3}(\Delta)}
\def\i#1{\int_{\lambda_1}^{\lambda_2}#1\,d\lambda}

\begin{document}
\title{Strichartz inequalities for the wave equation with the full
Laplacian on the
Heisenberg group}
\date{}

\author
{ Giulia FURIOLI \\
{\small Dipartimento di Ingegneria Gestionale e dell'Informazione, Universit\`a di
Bergamo,}\\
{\small Viale Marconi 5, I--24044 Dalmine (BG), Italy} \\
{\small E-mail:} \texttt{gfurioli@unibg.it},\\
\
\\
Camillo MELZI\\
{\small Dipartimento di Scienze Chimiche, Fisiche e Matematiche, Universit\`a
dell'Insubria,}\\
{\small Via Valleggio 11, I--22100 Como, Italy} \\
{\small E-mail:} \texttt{melzi@uninsubria.it}\\
\ \\
and Alessandro VENERUSO\\
{\small Dipartimento di Matematica, Universit\`a di Genova,}\\
{\small Via Dodecaneso 35, I--16146 Genova, Italy} \\
{\small E-mail:}  \texttt{veneruso@dima.unige.it}
}
\maketitle
\footnotetext{
{\em Keywords:} Strichartz inequalities, wave equation, full Laplacian, Heisenberg group.

{\em 2000 Mathematics Subject Classification}: 22E25, 35B65.}

\begin{abstract}
We prove dispersive and Strichartz inequalities for the solution of the wave
equation
related to the full
Laplacian on the Heisenberg group, by means of Besov spaces defined by a
Littlewood--Paley
decomposition related to the spectral resolution of the full Laplacian.
This requires a careful
analysis due also to the non-homogeneous nature of the full Laplacian.
This result has to be compared to a previous one  by Bahouri, G\'erard
and Xu concerning the solution of the wave equation related to
the Kohn-Laplacian.
\end{abstract}

\vfill
\eject
\section{Introduction}\leavevmode\par
The aim of this paper  is to study Strichartz inequalities for the
solution of the following Cauchy problem for the wave equation
on the Heisenberg group $\hn$ of topological dimension $2n+1$ and
homogeneous dimension $N=2n+2$:
\begin{equation}\label{cau}
\begin{cases}
        {\partial}^2_{ t}u  +\LL u= f\in L^1((0,T),
L^2(\H_n))& \\
        u(0)=u_0  \in \bpl 1 2 2&\\
        \partial_t u(0) =u_1 \in L^2(\H_n)
\end{cases}
\end{equation}
where $\LL$ is the full Laplacian
on $\hn$ (to be defined in Section
\ref{Notpre}) and the  Besov spaces $\bpl{\rho}{q}{r}$ are defined by a
Littlewood--Paley decomposition related to the spectral resolution of the
full Laplacian (see
Section \ref{LP-bes}). In
\cite{BGX}, Bahouri, G\'erard and Xu studied the analogous Cauchy
problem with the
Kohn-Laplacian $\Delta$ instead of the full Laplacian $\LL$,  using the
Besov spaces
$\bpd{\rho}{q}{r}$ which contain $\bpl{\rho}{q}{r}$ for $\rho>0$ (see
Proposition \ref{inclusioni}).
In
\cite{FV},  the first and last authors studied the corresponding Cauchy
problem for the
Schr\"odinger equation where they
introduced the full Laplacian instead of the Kohn-Laplacian, but still they
used the Besov spaces $\bpd{\rho}{q}{r}$.

Let us begin by recalling the structure of the solution of the Cauchy
problem \eqref{cau}.
It is well-known that the  solution of \eqref{cau} can be written as
$u=v+w$
where  $v$ is the  solution of \eqref{cau} with $f=0$ and $w$ is the
solution of \eqref{cau}
with $u_0=u_1=0$.
More precisely,
\begin{equation}
v(t)=  \cos (t\sqrt{\LL}) u_0+ \frac{\sin (t
\sqrt{\LL})}{\sqrt{\LL}}u_1\label{vxt}
\end{equation}
and
\begin{equation}
w(t)= \int_0^t \frac{\sin ((t-\sigma) \sqrt{\LL})}{\sqrt{\LL}} f(\sigma)\,
d\sigma.\label{wxt}
\end{equation}
We can now state the main results of this paper.
As always when dealing with Strichartz inequalities, we
prove first the following dispersive inequality on $v$.
\begin{prop}\label{Dispe}
Let $\rho\in[N-\frac{3}{2},N-\mez]$ and  $u_0\in \bpl{\rho}{1}{1}$, $
u_1\in \bpl
{\rho-1} 1 1$.
Then, there exists a constant
$C>0$, which does not depend on $u_0$, $u_1$, such that
\[
    \|v(t)\|_{L^{\infty}(\hn)} \leq C|t|^{-\mez}(\|u_0\|_{\bpl \rho 1
1}+\|u_1\|_{\bpl {\rho-1} 1 1
}),
    \qquad t\in\R^{*}.
\]
\end{prop}
Let us remark the main difference between
Proposition \ref{Dispe} and
\cite[Th\'eor\`eme 1.2]{BGX}: in the hypotheses of the latter theorem, they
obtain only the index $\rho=N-{\mez}$, which in
that case is sharp
because of the homogeneity property of the Kohn-Laplacian $\Delta$.

\medskip
For every interval $I\subset\R$ we will denote
by $L^p_I(X)$ the space $L^p(I,X)$. The Strichartz inequalities  we have
obtained are the
following.
\begin{theo}\label{str-in}
Let $r_1$, $r_2 \in [2, \infty]$. Let
$\rho_1$, $\rho_2 \in
\R$ and $p_1$,
$p_2\in [1, \infty]$
such that:
\begin{enumerate}[a)]
\item ${\dis  \frac 2{p_i} = \mez -  \frac 1{r_i}}$  for $i=1,2$;
\item $ -\left(N-\frac 12 \right )\left(\frac 12 - \frac 1{r_1}\right)+1 \leq
\rho_1 \leq -\left( N-\frac 32 \right)\left(\frac 12 - \frac 1{r_1}\right)+1$;
\item $ -\left(N-\frac 12 \right )\left(\frac 12 - \frac 1{r_2}\right) \leq
\rho_2 \leq -\left( N-\frac 32 \right)\left(\frac 12 - \frac 1{r_2}\right)$.
\end{enumerate}
Let $r'_i$, $p'_i$ such that ${\dis \frac 1{r'_i} + \frac 1{r_i} =1}$ and
${\dis \frac 1{p'_i} + \frac 1{p_i} =1}$ for $i=1,2$.
Then for every interval $I$ which contains $0$ the following estimates are
satisfied:
\begin{eqnarray*}
\|v\|_{\spc {p_1}{\R}{\rho_1}2{r_1}} +\|\partial_t v\|_{\spc
{p_1}{\R}{\rho_1-1}2{r_1}}
&\leq& C\,\left(\|u_0\|_{\bpl 1 2 2}+\|u_1\|_{\lp2}\right)\\
\| w\|_{\spc {p_1}{I}{\rho_1}2{r_1}} +  \| \partial_t w\|_{\spc
{p_1}{I}{\rho_1-1}2{r_1}}
&\leq& C\,\|f\|_{\spc {p'_2}{I}{-\rho_2}2{r'_2}}
\end{eqnarray*}
where the
constant $C>0$ depends neither on $u_0$, $u_1$, $f$ nor on the interval
$I$.
\end{theo}

So, we can  deduce from Theorem \ref{str-in} the following result, which we
compare to the
analogous result by Bahouri, G\'erard and Xu.
\begin{coro}\label{meglio}
Let $u$ be the solution of the Cauchy problem \eqref{cau}.
If $p$ and $r$ satisfy   $0\leq \frac 2p \leq \frac 12 -\frac 1r$ and
$(N-1)\left(\frac 12 -\frac 1r\right) -1 \leq \frac 1p \leq N\left( \frac
12 -\frac 1r\right) -1$, then there exists a
constant
$C \gt 0$, which does not depend on $u_0$, $u_1$, $f$, such that
for every interval $I$ which contains $0$ the following estimate is satisfied:
\[
\|u\|_{\sp {p}{I}{L^r(\hn)}}
\leq C\,\left(\|u_0\|_{\bpl 1 2 2}+\|u_1\|_{\lp2}
+\|f\|_{\sp{1}{I}{L^2(\hn)}}\right).
\]
\end{coro}
\begin{minipage}{7.5truecm}
In \cite[Th\'eor\`eme 1.1]{BGX}, the solution of the wave equation with the
Kohn-Laplacian was
proved to belong
to $\sp {p}{I}{L^r(\hn)}$ only for $p$ and $r$ satisfying $2N-1\leq p \leq
\infty$ and $\frac 1p =
N\left (\frac 12 -\frac 1r\right) -1$ which is a subset of the range of
values of $p$ and $r$ we
have found (since it is equivalent to
$0\leq \frac 2p \leq \frac 12 -\frac 1r$ and   $\frac 1p =  N\left (\frac
12 -\frac 1r\right) -1$).
The set of the admissible values  $\left(\frac 1r, \frac 1p\right)$ found
in Corollary \ref{meglio}
 is represented in the picture,
where the result by Bahouri, G\'erard and Xu corresponds to the segment  $BC$.

\end{minipage}
\ \
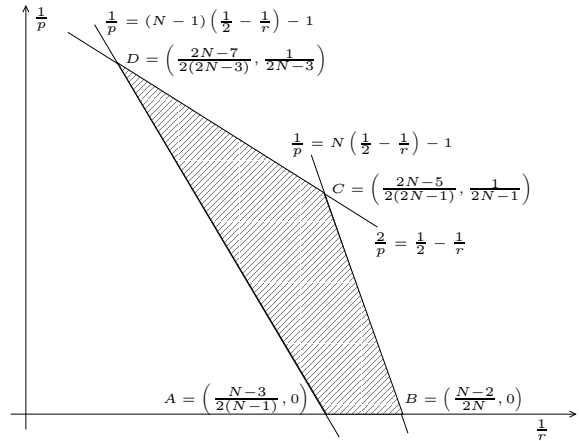
\begin{minipage}{6truecm}

\bigskip
\begin{picture}(0,0)%
\includegraphics{losanghina.pstex}%
\end{picture}%
\setlength{\unitlength}{1263sp}%
\begingroup\makeatletter\ifx\SetFigFont\undefined%
\gdef\SetFigFont#1#2#3#4#5{%
  \reset@font\fontsize{#1}{#2pt}%
  \fontfamily{#3}\fontseries{#4}\fontshape{#5}%
  \selectfont}%
\fi\endgroup%
\begin{picture}(11124,8499)(1789,-8773)
\put(8101,-3961){\makebox(0,0)[lb]{\smash{{\SetFigFont{5}{6.0}{\rmdefault}{\mddefault}{\updefault}{\color[rgb]{0,0,0}$C =\left(\frac {2N-5}{2(2N-1)}, \frac 1{2N-1}\right)$}%
}}}}
\put(2251,-586){\makebox(0,0)[lb]{\smash{{\SetFigFont{5}{6.0}{\rmdefault}{\mddefault}{\updefault}{\color[rgb]{0,0,0}$\frac 1p$}%
}}}}
\put(7276,-3061){\makebox(0,0)[lb]{\smash{{\SetFigFont{5}{6.0}{\rmdefault}{\mddefault}{\updefault}{\color[rgb]{0,0,0}$\frac 1p= N\left( \frac 12 -\frac 1r\right)-1$}%
}}}}
\put(4051,-1411){\makebox(0,0)[lb]{\smash{{\SetFigFont{5}{6.0}{\rmdefault}{\mddefault}{\updefault}{\color[rgb]{0,0,0}$D =\left(\frac {2N-7}{2(2N-3)}, \frac 1{2N-3}\right)$}%
}}}}
\put(9526,-8086){\makebox(0,0)[lb]{\smash{{\SetFigFont{5}{6.0}{\rmdefault}{\mddefault}{\updefault}{\color[rgb]{0,0,0}$B =\left(\frac {N-2}{2N}, 0\right)$}%
}}}}
\put(8910,-5016){\makebox(0,0)[lb]{\smash{{\SetFigFont{5}{6.0}{\rmdefault}{\mddefault}{\updefault}{\color[rgb]{0,0,0}$\frac 2p= \frac 12-\frac 1r$}%
}}}}
\put(3638,-666){\makebox(0,0)[lb]{\smash{{\SetFigFont{5}{6.0}{\rmdefault}{\mddefault}{\updefault}{\color[rgb]{0,0,0}$\frac 1p= (N-1)\left( \frac 12 -\frac 1r\right)-1$}%
}}}}
\put(4800,-8086){\makebox(0,0)[lb]{\smash{{\SetFigFont{5}{6.0}{\rmdefault}{\mddefault}{\updefault}{\color[rgb]{0,0,0}$A =\left(\frac {N-3}{2(N-1)},0\right)$}%
}}}}
\put(12076,-8686){\makebox(0,0)[lb]{\smash{{\SetFigFont{5}{6.0}{\rmdefault}{\mddefault}{\updefault}{\color[rgb]{0,0,0}$\frac 1r$}%
}}}}
\end{picture}%

\end{minipage}

\medskip
Other results on the sharpness of the dispersive inequalities and remarks
about the behaviour
of the operator $e^{-it\sqrt \LL}$ when analysed by the Besov spaces
$\bpd{\rho}{q}{r}$ can be found
in Section \ref{sharp}.
\section{Notation and preliminaries}\label{Notpre}\leavevmode\par
In this paper
$\N$ denotes the set of nonnegative integers, $\Z_{+}$ the set of positive
integers and
$\R_{+}$ the set of positive real numbers. For $p\in[1,\infty]$ we denote
by $p'$ the
conjugate index of $p$, such that
$\frac{1}{p}+\frac{1}{p'}=1$. We will denote by $C$ any positive
constant, depending only on the group, which will not be
necessarily the same at each occurrence.\par\smallskip
In this section we recall some basic facts about harmonic analysis on the
\gh. For the proofs and
further information, see e.g.
\cite{BJRW}, \cite{F}, \cite{G}, \cite{N}.\par\smallskip The \gh\ $\hn$,
$n\in\Z_{+}$, is the
nilpotent Lie group whose underlying \mbox{manifold} is
$\R^n\times\R^n\times\R$, with the following multiplication law:
\[(x,y,s)(x',y',s')=(x+x',y+y',s+s'+2(y\cdot x'-x\cdot y')),\quad
x,x',y,y'\in\R^n,\ s,s'\in\R.\]
The Lie algebra of $\hn$ is generated
by the left-invariant
vector fields
$\ppp{X}{n},\ppp{Y}{n},S$, where
\[ X_j=\ddd{x_j}+2{y_j}\ddd{s},\qquad
Y_j=\ddd{y_j}-2{x_j}\ddd{s},\qquad
S=\ddd{s}.\]
We indicate an element $g=(x,y,s)\in\hn$ as $g=(z,s)$, where
$z=x+iy\in\C^n$.  The family of dilations
$\{\delta_r:r>0\}$
given by
\[\delta_r(z,s)=(rz,r^2 s)\]
makes $\hn$ a stratified group of homogeneous dimension
$N=2n+2$. The Kohn-Laplacian
\[\Delta=-\sjn(X_j^2+Y_j^2)\]
satisfies the homogeneity property $\Delta(f\circ\delta_r)=r^2(\Delta
f\circ\delta_r)$, $r>0$, while
the full Laplacian
\[\LL=\Delta-S^2\]
is not invariant with respect to the dilation structure of $\hn$.\par\smallskip
The bi-invariant Haar measure $dg$ on $\hn$ coincides with the Lebesgue
measure on $\R^{2n+1}$. The convolution of two functions
$f_1$ and $f_2$ on $G$, defined by
\[f_1*f_2(g)=\ih f_1(g{g'}^{-1})f_2(g')\,dg',\quad g\in\hn,\]
satisfies the Young's inequality (where
$1+\frac{1}{r}=\frac{1}{p}+\frac{1}{q}$)
\[\|f_1*f_2\|_{\lp r}\leq\|f_1\|_{\lp p}\|f_2\|_{\lp q}.\]
The convolution of $\varphi\in\sh$ and $u\in\sih$, where $\sh$ is the Schwartz
space and $\sih$ is the space of tempered distributions, is defined as
usual (see e.g. \cite{V}).
We say that a function
$f$ on
$\hn$ is radial if
the value of $f(z,s)$
depends only on $|z|$ and $s$.
We denote by $\srad$ and by $\lrad p$,
$1\leq p\leq\infty$, the spaces of radial functions in $\sh$ and
in $\lp p$, respectively. The space $\lrad 1$ is a commutative, closed
$*$-subalgebra of $\lp 1$.
The Gelfand spectrum $\Sigma$ of $\lrad 1$
can be identified, as a measure space, with the space
$\N\times\R$ equipped with the Godement--Plancherel measure $\mu$
defined by
\[\int_{\Sigma} F(\psi)\,d\mu(\psi)=\idd{F\kl}.\]
The spherical Fourier
transform
of a function
$f\in\lrad 1$ is given by
\[\hat f\kl=\ih f(g)\okl(g)\,dg,\quad m\in\N,\ \lambda\in\R,\]
with
\[\okl(z,s)=\coefbi^{-1}e^{i\lambda s}e^{-|\lambda||z|^2}
\lmnu(2|\lambda||z|^2),\]
where $L_m^{(\alpha)}$ is the Laguerre polynomial of type $\alpha\in\N$ and
degree $m\in\N$,
defined by
\[L_m^{(\alpha)}(\tau)=\sum_{k=0}^m\frac{(-1)^k}{k!}\cb{m+\alpha}{k+\alpha}
\tau^k,\quad\tau\in\R.\] We have
$\widehat{f_1*f_2}=\widehat{f_1}\widehat{f_2}$ for any
$f_1,f_2\in\lrad 1$. Since
$\|\okl\|_{\lp\infty}=1$
the spherical Fourier transform is bounded from $\lrad 1$ to
$L^{\infty}(\Sigma)$.
Moreover,
by the Godement--Plancherel
theory, it extends uniquely to a unitary operator $\GG:\lrad
2\longrightarrow L^2(\Sigma)$. We still write $\hat f$ instead of $\GG
f$. If $f\in\lrad 2$ and
$\hat f\in L^1(\Sigma)$,
the following inversion formula holds:
\begin{equation}
f(g)=\idd{\hat
f\kl\omega_{m,-\lambda}(g)},\quad g\in\hn.\label{fx}
\end{equation}
The space $\GG(\srad)$ has been
described
in~\cite{BJR}. For our purposes, it is sufficient to remark that
$\GG(\srad)\subset L^1(\Sigma)$.
Moreover,
if $f\in\srad$ the functions $\Delta f$ and $\LL f$ are in $\srad$ and
their spherical Fourier
transforms are given by:
\begin{eqnarray}\widehat{\Delta f}\kl&=&4(2m+n)|\lambda|\hat
f\kl,\label{ddfml}\\
\widehat{\LL f}\kl&=&(4(2m+n)|\lambda|+\lambda^2)\hat
f\kl.\label{llfml}\end{eqnarray}
Both $\Delta$ and $\LL$ are positive self-adjoint operators densely defined
on $\lp 2$. So, by the
spectral theorem, for any bounded Borel function $h$ on
$[0,+\infty)$ the
operators $h(\Delta)$ and $h(\LL)$ are bounded on $\lp 2$. Since the point
0 may be
neglected in
the spectral resolution (see \cite{A}, \cite{C}), we consider
that
the function
$h$ is defined on $\R_{+}$. If $f\in\lrad 2$ the functions $h(\Delta)f$ and
$h(\LL)f$ are in $\lrad
2$ and their spherical Fourier transforms, by (\ref{ddfml}) and
(\ref{llfml}), are given by:
\begin{eqnarray}\widehat{h(\Delta)f}\kl&=&h(4(2m+n)|\lambda|)\hat
f\kl,\label{dfml}\\
\widehat{h(\LL)f}\kl&=&h(4(2m+n)|\lambda|+\lambda^2)\hat
f\kl.\label{lfml}\end{eqnarray}
If $f\in\srad$ then, by the previous remarks, the functions $h(\Delta)f$
and $h(\LL)f$ can be
recovered from their spherical Fourier transforms by means of the inversion
formula
(\ref{fx}).\par\smallskip The operators
$h(\Delta)$ and
$h(\LL)$ commute with left translations. So by the Schwartz' kernel
theorem, which is valid also on
$\hn$ (see \cite[Theorem 3.2]{KVW}),  they admit kernels in $\sih$, which
we call
$H_{\Delta}$ and
$H_{\LL}$ respectively, satisfying
$h(\Delta)f=f*H_{\Delta}$ and $h(\LL)f=f*H_{\LL}$ for any $f\in\sh$. If $h$
is the
restriction on
$\R_{+}$ of a function in ${\cal S}(\R)$, then $H_{\Delta}$ and $H_{\LL}$
are in $\srad$
(see \cite[Corollary 7]{FMV}; see also \cite{H}, \cite{M} for $H_{\Delta}$,
\cite{V} for $H_{\LL}$)
and their spherical Fourier transforms, by
(\ref{dfml}) and
(\ref{lfml}), are
given by:
\begin{eqnarray*}\widehat{H_{\Delta}}\kl&=&h(4(2m+n)|\lambda|),\\
\widehat{H_{\LL}}\kl&=&h(4(2m+n)|\lambda|+\lambda^2).\end{eqnarray*}

\section{Littlewood--Paley decompositions and Besov spaces}
\label{LP-bes}\leavevmode\par
Let $R$ be a non-negative function in $C^{\infty}(\R)$ such
that supp$\,R\subset[\frac{1}{4},4]$ and
\[\sum_{j\in\Z}R(2^{-2j}\tau)=1,   \quad \tau>0.\]
For any $j\in\Z$ we denote by $\fj$ and $\psi_j$ the kernels of the
operators $R(\dj\Delta)$ and
$R(\dj\LL)$, respectively.
The remarks at the end of Section~\ref{Notpre}
guarantee that
$\fj,\psi_j\in\srad$ and
\begin{eqnarray}\widehat{\fj}\kl&=&R(2^{2-2j}(2m+n)|\lambda|),\label{fjml}\\
\widehat{\psi_j}\kl&=&R(\dj(4(2m+n)|\lambda|+\lambda^2)).\label{psjml}\end{eqnarray}
If $j,k\in\Z$ with $|j-k|\geq 2$, then
$\varphi_j*\varphi_{k}=\psi_j*\psi_{k}=0.$ Moreover
we have the following
\begin{lem}
For any $j\in\Z$ the sets $U_j=\{k\in\Z:\varphi_j*\psi_k\neq 0\}$
and
$V_j=\{k\in\Z:\psi_j*\varphi_k\neq 0\}$ are finite and
$\min U_j\geq j-2$, $\max V_j\leq j+2.$\label{lemmaprimo}\end{lem}
\begin{proof} Fix $j\in\Z$ and
$k\in U_j$. By (\ref{fjml}) and (\ref{psjml}) there exist $m\in\N$,
$\lambda\in\R$ such that
\[R(2^{2-2j}(2m+n)|\lambda|)R(2^{-2k}(4(2m+n)|\lambda|+\lambda^2))\neq
0.\] Put $\xi=4(2m+n)|\lambda|$ and $\eta=\lambda^2$. The pair
$(\xi,\eta)$ satisfies the following system of inequalities:
\begin{equation}\left\{\begin{array}{ll}\frac{1}{4}\leq 2^{-2j}\xi\leq 4\\
\frac{1}{4}\leq 2^{-2k}(\xi+\eta)\leq4\\
0\leq\eta\leq\frac{\xi^2}{16n^2}.\end{array}\right.\label{system}\end{equation}
On the other hand, it is easy to check that the system (\ref{system}) admits
solutions only if
\[2^{2j-4}\leq 2^{2k}\leq\frac{2^{4j+2}}{n^2}+2^{2j+4}.\]
These conditions give the conclusion not only for $U_j$, but also for
$V_j$: for the latter one it
is sufficient to interchange the roles of $j$ and $k$, noting that $k\in
V_j$ if and only if
$j\in U_k$.\end{proof}\par\medskip
A direct application of the inversion
formula~(\ref{fx}) gives
\begin{equation}\fj(z,s)= 2^{Nj} \varphi_0(2^{j}z, 2^{2j}s),\quad j\in\Z,\
(z,s)\in\hn.\label{homprop}\end{equation} So
\begin{equation} \|\fj\|_{\lp 1}=\|\varphi_0\|_{\lp 1},\quad
j\in\Z.\label{hom}\end{equation}
On the other hand, despite the lack of homogeneity, by \cite[Proposition
6]{FMV}
there exists $C>0$ such
that
\begin{equation} \|\psi_j\|_{\lp 1}\leq C,\quad
j\in\Z.\label{normauno}\end{equation}
In this section, in order to carry on some results which
are valid for both operators $\Delta$ and $\LL$, we use the notation $L$ to
denote either
$\Delta$ or $\LL$. For any $u\in\sih$, if
$L=\Delta$ we set
$\Delta_j u=u*\fj$, if $L=\LL$ we set $\Delta_j u=u*\psi_j$. By standard
arguments (see e.g.
\cite[Proposition 9]{FMV}) we can deduce from (\ref{hom}) and
(\ref{normauno}) that
\begin{equation}\|L^{\frac{\sigma}{2}}\Delta_j u\|_{\lp p}\leq C
2^{j\sigma}\|\Delta_j u\|_{\lp p},\quad \sigma\in\R,\ j\in\Z,\
1\leq p\leq\infty,\ u\in\sih,\label{fmv}\end{equation} where both sides of
(\ref{fmv}) are allowed
to be infinite.\par\medskip By the spectral theorem, for
any
$f\in\lp 2$ the following homogeneous
\LP decomposition holds:
\begin{equation}f=\sjz \Delta_j f\quad{\rm in}\ \lp
2.\label{sdjf}\end{equation}
So
\begin{equation}\|f\|_{\lp{\infty}}\leq\sjz\|\Delta_j f\|_{\lp{\infty}},\quad
f\in\lp
2,\label{fv}\end{equation} where both sides of (\ref{fv}) are allowed to be
infinite.\par\medskip
The methods of
\cite{St}, together with any multiplier theorem for $L$ (see \cite{A}; see
also \cite{He},
\cite{MS} for $L=\Delta$, \cite{MRS1}, \cite{MRS2} for $L=\LL$), yield the
following \LP theorem:
\begin{prop}\label{litpal}
Let $1<p<\infty$ and $u\in\sih$. The following facts are equivalent:
\begin{enumerate}[(i)]
\item $u\in\lp p$;
\item $u=\sjz \Delta_j u$ in $\sih$ and $(\sjz|\Delta_j u|^2)^{\mez}\in\lp p$.
\end{enumerate}
Moreover, if $u\in\lp p$ then
\[\|u\|_{\lp
p}\sim\|(\sjz|\Delta_j u|^2)^{\mez}\|_{\lp p}.\]
\end{prop}
{\em Remark:}\quad For $L=\Delta$ Proposition \ref{litpal} has been proved
also in \cite[Proposition
2.3]{BGX} using the homogeneity property (\ref{homprop}).\par\medskip
Let $q,r\in [1,\infty]$ and $\rho\in\R$;
the homogeneous Besov space $\bp{\rho}{q}{r}$ associated to the operator
$L$ is defined as follows:
\[\bp{\rho}{q}{r}=\{u\in \sih: u=\sjz \Delta_j u\ {\rm in\ }\sih\ {\rm and\ }
\{2^{j\rho}\|\Delta_j u\|_{\lp r} \}_{j\in\Z}\in l^q(\Z)\}.\]
We collect in the following proposition all the properties we need about
the spaces
$\bp{\rho}{q}{r}$.
\begin{prop}\label{besovprop} Let $q,r\in [1,\infty]$ and $\rho<\frac{N}{r}$.
\begin{enumerate}[(i)]
\item The space $\bp{\rho}{q}{r}$ is a Banach space endowed with the norm
\[\|u\|_{\bp{\rho}{q}{r}} =\left\|\{2^{j\rho} \|\Delta_j u\|_{\lp
r}\}_{j\in\Z}\right\|_{l^q(\Z)};\]
\item the definition of $\bp{\rho}{q}{r}$ does not depend on the choice of
the function $R$ in the \LP decomposition;
\item for any $u\in\sih$ and $\sigma>0$ we have that $u\in\bp{\rho}{q}{r}$
if and only if
$L^{\frac{\sigma}{2}} u\in\bp{\rho-\sigma}{q}{r}$, with
\[\|u\|_{\bp{\rho}{q}{r}}\sim\|L^{\frac{\sigma}{2}}u\|_{\bp{\rho-\sigma}{q}{r}};
\]
\item the inclusion $\bp{\rho}{q}{r}\subset\sih$ is continuous;
\item if $-\frac{N}{r'}<\rho<\frac{N}{r}$ then $\sh\subset\bp{\rho}{q}{r}$ with
continuous inclusion;
\item if $q,r\in [1,\infty)$ and $-\frac{N}{r'}<\rho<\frac{N}{r}$, then $\sh$
is dense in $\bp{\rho}{q}{r}$;
\item if $q,r\in [1,\infty)$ and $-\frac{N}{r'}<\rho<\frac{N}{r}$, the dual
space
of
$\bp{\rho}{q}{r}$
is $\bp{-\rho}{q'}{r'}$;
\item for all $q\in [1, \infty]$ and $\alpha \in [N-1,N]$ we have the
continuos inclusions
\[\bpl {\rho_1}{q}{r_1} \subset \bpl {\rho_2}q{r_2}, \quad {\dis \frac 1{r_1} -
\frac {\rho_1}{\alpha} =
\frac 1{r_2} - \frac {\rho_2}{\alpha} }, \ \rho_1\geq \rho_2;\]
\[\bpd {\rho_1}{q}{r_1} \subset \bpd {\rho_2}q{r_2}, \quad {\dis \frac 1{r_1} -
\frac {\rho_1}{N} =
\frac 1{r_2} - \frac {\rho_2}{N} }, \ \rho_1\geq \rho_2;\]
\item for all $r\in [2, \infty)$ we have the continuos inclusion
$\bp 0 2 r\subset\lp r$;
\item $\bp 022 = \lp 2$ with equivalent norms;
\item for all\quad  $\vartheta, \rho_1, \rho_2, q_1, q_2, r_1, r_2$  satisfying
${\displaystyle \vartheta \in [0,1], \ q_i, r_i\in(1,\infty),\
\rho_i<\frac{N}{r_i}}$,
we have
\[
[\bp {\rho_1}{q_1}{r_1}, \bp {\rho_2}{q_2}{r_2}]_\vartheta  = \bp {\rho}{q}{r}
\]
with ${\displaystyle \rho = (1-\vartheta) \rho_1 +\vartheta \rho_2}$,\
 ${\displaystyle \frac 1q = \frac {1-\vartheta}{q_1}+\frac
{\vartheta}{q_2}}$ and\
 ${\displaystyle \frac 1r = \frac {1-\vartheta}{r_1}+
\frac {\vartheta}{r_2}}$.
\end{enumerate}
\end{prop}
We omit the proof of Proposition \ref{besovprop}. In fact, all the
statements of the proposition are well-known for the spaces $\bpd \rho q r$
(see \cite{BG}, \cite{BGX},
\cite{FV}) and the proofs for the spaces $\bpl \rho q r$ are analogous: the
only properties really
needed are estimates (\ref{normauno}) and (\ref{fmv}), Proposition
\ref{litpal} and the fact that the kernel of $h(\LL)$ is in $\sh$ if $h\in
{\cal S}(\R)$ (see Section~\ref{Notpre}). Once we
have these properties, we can prove Proposition~\ref{besovprop} by the
methods in
\cite{P}, which do not involve any homogeneity property.
More generally,
we could define
homogeneous Besov spaces and prove, with the same methods, an analogous proposition
 in the more general context of a
nilpotent Lie group $G$
endowed with a
sub-Laplacian
$L=-\sum_{j=1}^k X^2_j$, where
$\ppp{X}{k}$ are left-invariant vector fields on $G$ which satisfy
the H\"ormander's condition, i.e. they
generate, together with their successive Lie brackets
$[X_{i_1},[\ldots,X_{i_{\alpha}}]\cdots]$,
the Lie algebra of $G$.
For more details about properties of Besov spaces
in this context, see
\cite{S1}, \cite{S2}, \cite{FMV}, where nevertheless inhomogeneous Besov
spaces are considered.
Here we want to prove some continuous inclusions
between the two kinds of
homogeneous Besov spaces which we have introduced.
\begin{prop} The following continuos inclusions hold:
\begin{eqnarray}\bpl \rho q r&\subset&\bpd \rho q r,\quad 1\leq
q\leq\infty,\ 1\leq
r<\infty,\
0<\rho<\frac{N}{r};\label{incluno}\\
\bpd \rho q r&\subset&\bpl \rho q r,\quad 1\leq q\leq\infty,\ 1<r\leq\infty,\
-\frac{N}{r'}<\rho<0.\label{incldue}\end{eqnarray}\label{inclusioni}\end{prop}
\begin{proof} We only prove (\ref{incluno}), since the proof of
(\ref{incldue}) is analogous. Fix
$u\in\bpl \rho q r$, with $1\leq q\leq\infty,\ 1\leq
r<\infty$ and
$0<\rho<\frac{N}{r}$. Since $u=\sjz u*\psi_j$ in $\sih$, by Lemma
\ref{lemmaprimo} we have
$u*\varphi_k=\sum_{j\geq k-2}
u*\psi_j*\varphi_k$ in
$\sih$ for any
$k\in\Z$, and so
\[
\begin{aligned}
2^{k\rho}\|u*\varphi_k\|_{\lp r} &\leq 2^{k\rho} \sum_{j\geq k-2}
\|u*\psi_j*\varphi_k\|_{\lp r}\\
&\leq C \sum_{j\geq k-2}2^{(k-j)\rho}2^{j\rho} \|u*\psi_j\|_{\lp r}
\end{aligned}
\]
by (\ref{hom}).
Therefore,  by Young's inequality
\[
\left\|\{2^{k\rho} \|u*\varphi_k\|_{\lp
r}\}_{k\in\Z}\right\|_{l^q(\Z)}\leq C \|u\|_{\bpl \rho q r}.
\]
We still have to prove that $u=\sum_{k\in\Z}u*\varphi_k$ in $\sih$.
By
Lemma
\ref{lemmaprimo}, for any
$f\in\sh$ we have:
\beq\sjz\sum_{k\in\Z}|\pair{u*\psi_j*\varphi_k}{f}|&=
&\sjz\sum_{h=-\infty}^2|\pair{u*\psi_j}{f*\varphi_{j+h}}|\\
&\leq&\sum_{h=-\infty}^2 2^{h\rho}\left(\sjz 2^{j\rho}\|u*\psi_j\|_{\lp
r}2^{-(j+h)\rho}\|f*\varphi_{j+h}\|_{\lp{r'}}\right)\\
&\leq&\sum_{h=-\infty}^2 2^{h\rho}\|u\|_{\bpl \rho q
r}\|f\|_{\bpd{-\rho}{q'}{r'}}<+\infty.\eeq
Note that $\varphi_k=\sjz\varphi_k*\psi_j$ in $\sh$ for any
$k\in\Z$, by (\ref{sdjf}) and Lemma \ref{lemmaprimo}.
Therefore, since $u=\sjz u*\psi_j$ in $\sih$, by Fubini's theorem we have
\[\pair{u}{f}=\sjz\sum_{k\in\Z}\pair{u*\psi_j*\varphi_k}{f}=\sum_{k\in\Z}\pair{u
*\varphi_k}{f},\quad
f\in\sh.\]
\end{proof}
However,
with the exception of particular cases as
$\rho=0$, $q=r=2$ (see Proposition
\ref{besovprop} (x)), the spaces $\bpd \rho q r$ and $\bpl \rho q r$ do not
coincide: for example, by
applying the Godement--Plancherel's formula and arguing as in the proof of
Lemma~\ref{lemmaprimo},
it is not hard to check that for $j\rightarrow
+\infty$ we have
\[\|\fj\|_{\bpd \rho q 2}\sim 2^{j(\rho+\frac{N}{2})},\qquad \|\fj\|_{\bpl
\rho q 2}\sim
2^{j(2\rho+\frac N 2)},\qquad 1\leq q\leq\infty,\ 0<\rho<\frac{N}{2}.\]
As a further evidence, in the following we will see that the spaces
$\bpd \rho q r$ and $\bpl \rho q
r$ have a very different behaviour with respect to Strichartz estimates
for the
solution of the Cauchy problem (\ref{cau}).
\section{Dispersive estimates}\label{Disp}\leavevmode\par
We begin by proving Proposition \ref{Dispe}. Let us introduce the tools of
the method; first of
all, we recall
the stationary phase lemma (see e.g. \cite{St}, pages 332--334) that will
be the central argument:

\begin{lem} Suppose $g,h\in C^{\infty}([a,b])$, with $g$ real-valued and
$h(b)=0$.
Suppose also
$|g^{(k)}(x)|\geq\delta$ for any
$x\in [a,b]$, with
$k\in\Z_{+}$ and $\delta>0$. If
$k=1$, we also require  that $g'$ is monotonic in $[a,b]$. Then there exists
a constant $C_k>0$, which
depends only on $k$ but not on $a,b,g,h,\delta$, such that
\[\left|\int_a^b e^{-ig(x)}h(x)\,dx\right|\leq
C_k\delta^{-\frac{1}{k}}\int_a^b|h'(x)|\,dx.\]
\label{ondeuno}\end{lem}

Moreover, we will use the following properties of the Laguerre polynomials
(see \cite{BGX}, \cite{EMOT}):

\begin{lem} Fix $\alpha\in\N$. There exists $C_{\alpha}>0$ such that for
$\tau\geq 0$ and
$m\in\N$ we have:
\[\left|L_m^{(\alpha)}(\tau)e^{-\frac{\tau}{2}}\right|\leq
C_{\alpha}(m+1)^{\alpha},\qquad
\left|\tau\frac{d}{d\tau}\left(L_m^{(\alpha)}(\tau)
e^{-\frac{\tau}{2}}\right)\right|\leq
C_{\alpha}(m+1)^{\alpha}.\]
\label{ondedue}\end{lem}

Finally, we will exploit the following estimates, which can be easily proved by
comparing the sums with
the corresponding integrals:

\begin{lem} Fix $\beta\in\R$. There exists
$C_{\beta}>0$
such that for $0<a<b$ and $n\in\Z_{+}$ we have:
\begin{eqnarray}
\sum_{\stackrel{\scriptstyle m\in\N}{2m+n\geq
a}}(2m+n)^{\beta}&\leq&C_{\beta}a^{\beta+1},\quad\beta<-1;\label{ein}\\
\sum_{\stackrel{\scriptstyle m\in\N}{2m+n\leq
b}}(2m+n)^{\beta}&\leq&C_{\beta}b^{\beta+1},\quad\beta>-1;\label{zwei}\\
\sum_{\stackrel{\scriptstyle m\in\N}{a\leq 2m+n\leq
b}}(2m+n)^{-1}&\leq&\log(C\frac{b}{a}).\label{trei}
\end{eqnarray}
\label{ondequattro}\end{lem}

We can now prove the following

\begin{prop}\label{ondecinque}
There exists a constant $C>0$, which depends only on $n$,
such that for any $\rho\in[N-\frac{3}{2},N-\mez]$, $j\in\Z$ and
$t\in\R^{*}$ we have:
\[\|e^{-it\sqrt{\LL}}\psi_j\|_{\lp \infty}\leq
C|t|^{-\mez}2^{j\rho}.\]
\end{prop}
\begin{proof}
Fix $t\in\R^{*}$, $j\in\Z$ and $\zs\in\hn$. By (\ref{fx}), (\ref{lfml}) and
(\ref{psjml}), putting $\sigma=\frac{s}{t}$ and $M=2m+n$ inside the sum
over $m$, we have
\beq\lefteqn{e^{-it\sqrt{\LL}}\psi_j\zs}\\
&=&\cn\smo\ir{e^{-it(\sigma\lambda
+\sqrt{4M|\lambda|+\lambda^2})}R(\dj(4M|\lambda|+\lambda^2))
e^{-|\lambda||z|^2}\lmnu(2|\lambda||z|^2)}.\eeq
Performing the change of variable
$x=\dj M\lambda$ we obtain
\[e^{-it\sqrt{\LL}}\psi_j\zs=\cn\dnj\smo\int_{\R}e^{-it\ddj
g_{j,\sigma,m}(x)}h_{j,z,m}(x)\,dx\]
where
\begin{eqnarray}g_{j,\sigma,m}(x)&=&\frac{1}{M}\left(\sigma
x+\sqrt{2^{2-2j}M^2|x|+x^2}\right),\label{gjsx}\\
h_{j,z,m}(x)&=&R\left(4|x|+\frac{\ddj
x^2}{M^2}\right)e^{-\frac{\ddj|x||z|^2}{M}}
\lmnu\left(\frac{2^{1+2j}|x||z|^2}{M}\right)\frac{|x|^n}{M^{n+1}}.
\label{hzmx}\end{eqnarray}
So
\begin{equation}{\rm supp}\,h_{j,z,m}\subset\{x\in\R:\frac{1}{4}\leq
4|x|+\frac{\ddj x^2}{M^2}\leq
4\}=\{x\in\R:a_{j,m}\leq|x|\leq b_{j,m}\}\label{supphzm}\end{equation}
where
\[a_{j,m}=\frac{1}{8(1+\sqrt{1+2^{2j-4}M^{-2}})},\qquad
b_{j,m}=\frac{2}{1+\sqrt{1+2^{2j}M^{-2}}}.\]
In particular
\begin{equation}b_{j,m}\leq\min\{1,2^{1-j}M\}.\label{bm}\end{equation}
Note that $g_{j,\sigma,m}(-x)=g_{j,-\sigma,m}(x)$ and
$h_{j,z,m}(-x)=h_{j,z,m}(x)$. Therefore, by symmetry we can consider only
the integrals
\[I_m=\int_{a_m}^{b_m} e^{-it\ddj
g_m(x)}h_m(x)\,dx\]
where we write $g_m,h_m,a_m,b_m$
for
$g_{j,\sigma,m},h_{j,z,m},a_{j,m},b_{j,m}$ respectively.
We prove that
\begin{equation}\smo|I_m|\leq\left\{\begin{array}{ll}C|t|^{-\mez}2^{-\frac{3}{2}
j},\quad &j\geq 0,\\
C|t|^{-\mez}2^{-\frac{j}{2}},\quad
&j<0.\end{array}\right.\label{sim}\end{equation}
For $x\in[a_m,b_m]$, by (\ref{gjsx}) we have
\begin{eqnarray}g'_m(x)&=&\frac{1}{M}\left(\sigma+\sqrt{1+\frac{2^{2-4j}M^4}{2^{
2-2j}M^2
x+x^2}}\right),\label{gim}\\
g''_m(x)&=&-2^{2-4j}M^3(2^{2-2j}M^2
x+x^2)^{-\met 3}.\label{giim}\end{eqnarray}
Note that by (\ref{supphzm}) we have
\begin{equation}2^{-2-2j}M^2\leq 2^{2-2j}M^2
x+x^2\leq 2^{2-2j}M^2,\quad x\in[a_m,b_m].\label{facile}\end{equation}
So (\ref{giim}) and (\ref{facile}) yield
\begin{equation}2^{-1-j}\leq|g''_m(x)|\leq 2^{5-j},\quad
x\in[a_m,b_m].\label{baasta}\end{equation}
Furthermore,
by Lemma~\ref{ondedue} and (\ref{bm}), one can verify that
\begin{equation}\|h'_m\|_{L^1([a_m,b_m])}\leq
\left\{\begin{array}{ll}C2^{-nj}M^{n-2},\quad &M\leq 2^j,\\
CM^{-2},\quad &M>2^j.\end{array}\right.\label{ihp}\end{equation}
So, by Lemma~\ref{ondeuno} with $k=2$, we obtain
\begin{equation}|I_m|\leq
\left\{\begin{array}{ll}C\tm 2^{-(n+\frac{1}{2})j}M^{n-2},\quad
&M\leq 2^j,\\
C\tm 2^{-\met j}M^{-2},\quad
&M>2^j.\end{array}\right.\label{stimaim}\end{equation}
For $j<0$, (\ref{sim}) follows directly from (\ref{stimaim}). For $n\geq 2$
and $j\geq 0$,
(\ref{sim}) still follows from (\ref{stimaim}) by applying
Lemma~\ref{ondequattro}
separately to the sums $\sum_{M\leq 2^j}|I_m|$ and $\sum_{M>2^j}|I_m|$. But
for
$n=1$ and $j\geq 0$ this argument does not work, since we cannot apply
(\ref{zwei}) to the sum
$\sum_{M\leq 2^j}|I_m|$.\par\smallskip
So
from now on we assume
$n=1$ and $j\geq 0$. We divide $\N$ into five (possibly empty) disjoint
subsets:
\beq
A_1&=&\{m\in\N:M>2^j\},\\
A_2&=&\{m\in\N:M\leq 2^j,M\leq\tm 2^{\met j}\},\\
A_3&=&\{m\in\N:M\leq 2^j,M>\tm 2^{\met j},\sigma\geq-\sqrt{1+2^{-1-2j}M^2}\},\\
A_4&=&\{m\in\N:M\leq 2^j,M>\tm 2^{\met j},\sigma\leq-\sqrt{1+2^{5-2j}M^2}\},\\
A_5&=&\{m\in\N:M\leq 2^j,M>\tm 2^{\met
j},-\sqrt{1+2^{5-2j}M^2}<\sigma<-\sqrt{1+2^{-1-2j}M^2}\}.\eeq
Then our assertion reads:
\begin{equation}\sum_{m\in A_r}|I_m|\leq C|t|^{-\mez}2^{-\frac{3}{2}j},\quad
r=1,\ldots,5,\ A_r\neq\emptyset.\label{smar}\end{equation}
We prove (\ref{smar}) separately for each $r$, using
each time Lemma
\ref{ondequattro}: precisely, we will use (\ref{ein}) for $r=1,3,4$,
(\ref{zwei}) for $r=2$ and
(\ref{trei}) for $r=5$. The case
$r=1$ can be treated as for
$n\geq 2$. For
$r=2$ we estimate
$\sum_{M\leq\tm 2^{\met j}}|I_m|$ by means of the inequality
\[|I_m|\leq C\|h_m\|_{L^1([a_m,b_m])}\leq C\dj\]
which follows from (\ref{hzmx}), (\ref{bm}) and Lemma \ref{ondedue}.
For $r=3,4$ we estimate $\sum_{M>\tm 2^{\met
j}}|I_m|$ by means of Lemma \ref{ondeuno} applied with $k=1$, using
(\ref{ihp}) and the estimates
\beq g'_m(x)&\geq&\frac{1}{M}\left(\sqrt{1+\dj
M^2}-\sqrt{1+2^{-1-2j}M^2}\right)\geq C\dj
M,\quad m\in A_3,\\
-g'_m(x)&\geq&\frac{1}{M}\left(\sqrt{1+2^{5-2j}M^2}-\sqrt{1+2^{4-2j}M^2}\right)\
geq C\dj
M,\quad m\in A_4,\eeq
which are consequences of (\ref{gim}) and (\ref{facile}). For $r=5$ we note
that
$A_5\neq\emptyset$ implies $\sigma<-1$ and $M\in J=(2^{j-\met
5}\sqrt{\sigma^2-1},2^{j+\met 1}\sqrt{\sigma^2-1})$. Then we estimate
$\sum_{M\in J}|I_m|$ by means
of Lemma \ref{ondeuno} applied with $k=2$.
\end{proof}

From Proposition \ref{ondecinque} we can obtain, by the same proof as in
\cite[pages 114--115]{BGX}, \cite[Corollary 10]{FV},
the following
\begin{coro} For $\rho\in[N-\frac{3}{2},N-\mez]$ there exists a constant
$C_{\rho}>0$ such that
\begin{align}
\left\|e^{-it\sqrt{\LL}}f \right\|&_{\lp \infty}\leq C_{\rho}
|t|^{-\mez}\|f\|_{\bpl{\rho}{1}{1}},\qquad f\in\sh,\ t\in\R^{*},\label{sopra}\\
\left\|e^{-it\sqrt{\LL}}f \right\|&_{\bpl {-1} 1 \infty}\leq C_{\rho}
|t|^{-\mez}\|f\|_{\bpl{\rho-1}{1}{1}},\qquad f\in\sh,\
t\in\R^{*}.\label{sotto}\end{align}
\label{ondesei}\end{coro}

The proof of  the dispersive inequality  is now straightforward.

\medskip
{\noindent{\bf Proof of Proposition  \ref{Dispe}:} }
By \eqref{sopra} we obtain
\begin{equation*}
    \|\cos t\sqrt{\LL} u_0\|_{L^{\infty}(\hn)} \leq
C|t|^{-\mez}\|u_0\|_{\bpl \rho 1 1}
\end{equation*}
and by \eqref{fv}, \eqref{fmv} and \eqref{sotto} we obtain
\begin{align*}
    \left\|\frac{\sin t \sqrt{\LL}}{\sqrt{\LL}}u_1\right\|_{L^\infty(\hn)}
       & \leq \sum_{j \in \Z}\left\|\frac{\sin t
\sqrt{\LL}}{\sqrt{\LL}}u_1\ast\psi_j\right\|_{L^\infty(\hn)}  \\
       & \leq \sum_{j \in \Z}2^{-j}\|\sin t \sqrt{\LL}
u_1\ast\psi_j\|_{L^\infty(\hn)}
         \leq C|t|^{-\mez}\|u_1\|_{\bpl {\rho-1} 1 1 }.
\end{align*}
\hfill $\scriptstyle\bullet$\vskip\theorempostskipamount\par

\section{Strichartz inequalities}\label{str}\leavevmode\par
We can now prove  Theorem \ref{str-in} and Corollary \ref{meglio}.

{\noindent{\bf Proof of Theorem \ref{str-in}: }}

By \eqref{vxt}  we can write
\[
\partial_t v(t)=   -\frac{e^{it\sqrt{\LL}}-e^{-it\sqrt{\LL}}}{2i}  \sqrt{\LL}
u_0+\frac{e^{it\sqrt{\LL}}+e^{-it\sqrt{\LL}}}{2}u_1
\]
where $\sqrt{\LL} u_0$ and $u_1$ both belong to $L^2(\hn)$.
Analogously, by \eqref{wxt}
\[
\partial_t w(t)
=\int_0^t\frac{e^{i(t-\sigma)\sqrt{\LL}}+e^{-i(t-\sigma)\sqrt{\LL}}}{2}
f(\sigma)\, d\sigma.
\]
So
\begin{equation*}
\begin{split}
\ &\|v\|_{\spc {p_1}{\R}{\rho_1}2{r_1}} +\|\partial_t v\|_{\spc
{p_1}{\R}{\rho_1-1}2{r_1}} \\
\leq &\, C\,
\left(\|\sqrt{\LL} \, v\|_{\spc {p_1}{\R}{\rho_1-1}2{r_1}} +\|\partial_t
v\|_{\spc
{p_1}{\R}{\rho_1-1}2{r_1}}\right)\\
 \leq &\, C\,  \left(\| e^{-it\sqrt{\LL}} \sqrt{\LL} u_0\|_{\spc
{p_1}{\R}{\rho_1-1}2{r_1}}+
\|e^{-it\sqrt{\LL}}u_1\|_{\spc {p_1}{\R}{\rho_1-1}2{r_1}}\right)
\end{split}
\end{equation*}
and
\begin{equation*}
\begin{split}
&\| w\|_{\spc {p_1}{I}{\rho_1}2{r_1}} +  \| \partial_t w\|_{\spc
{p_1}{I}{\rho_1-1}2{r_1}} \\
& \leq C\,\left(
\| \sqrt{\LL}\, w\|_{\spc {p_1}{I}{\rho_1-1}2{r_1}} +  \| \partial_t w\|_{\spc
{p_1}{I}{\rho_1-1}2{r_1}}\right)\\
&  \leq C\, \left(\left\| \int_0^t e^{i(t-\sigma) \sqrt{\LL}}
f(\sigma)\, d\sigma
\right\|_{\spc {p_1}{I}{\rho_1-1}2{r_1}}+\left\| \int_0^t e^{-i(t-\sigma)
\sqrt{\LL}}
f(\sigma)\, d\sigma
\right\|_{\spc {p_1}{I}{\rho_1-1}2{r_1}}\right).
\end{split}
\end{equation*}
Theorem \ref{str-in} follows therefore easily  by the following one, where
we have renamed $\rho_1$
the value $\rho_1-1$.
\begin{theo}\label{str-sempl}
Let $r_1$, $r_2 \in [2, \infty]$. Let
$\rho_1$, $\rho_2 \in
\R$ and $p_1$,
$p_2\in [1, \infty]$
such that:
\begin{enumerate}[a)]
\item ${\dis  \frac 2{p_i} = \mez -  \frac 1{r_i}}$  for $i=1,2$;
\item $ -\left(N-\frac 12 \right )\left(\frac 12 - \frac 1{r_i}\right) \leq
\rho_i \leq -\left( N-\frac 32 \right)\left(\frac 12 - \frac 1{r_i}\right)$
for $i=1,2$.
\end{enumerate}
Let $r'_i$, $p'_i$ such that ${\dis \frac 1{r'_i} + \frac 1{r_i} =1}$ and
${\dis \frac
1{p'_i} + \frac 1{p_i} =1}$ for $i=1,2$.
Then for every interval $I$ which contains $0$ the following estimates are
satisfied:
\begin{eqnarray*}
\|e^{ -it\sqrt{\LL}} u_0\|_{\spc
{p_1}{\R}{\rho_1}2{r_1}}&\leq&C\,\|u_0\|_{\lp2}\\
\left \| \int_0^t e^{\pm i
(t-\sigma)\sqrt{\LL}} f(\sigma)\,d\sigma\right \|_{\spc {p_1}{
I}{\rho_1}2{r_1}}
&\leq&C\,\|f\|_{\spc {p'_2}{I}{-\rho_2}2{r'_2}}
\end{eqnarray*}
where the
constant $C>0$ depends neither on $u_0$, $f$ nor on the interval $I$.
\hfill $\scriptstyle\bullet$\vskip\theorempostskipamount\par
\end{theo}

We omit the proof of Theorem~\ref{str-sempl}: in fact, once we have obtained
Proposition~\ref{ondecinque}, the procedure is classical and a good
reference is given, for
example, by the  papers by Ginibre and Velo
(\cite{GV2}) or by Ginibre (\cite{Gi}).
A detailed presentation in this framework is also given by \cite{FV}.

{\noindent{\bf Proof of Corollary \ref{meglio}: }}
Let us remark first that for $\rho_1 \geq 0$, Proposition \ref{besovprop}
((viii) and (ix)) implies
\begin{equation}\label{sob}
\bpl{\rho_1}2{r_1} \subset \bpl {0}2{r_{\min}} \cap  \bpl
{0}2{r_{\max}}\subset\lp{r_{\min}}\cap\lp{r_{\max}},
\end{equation}
where $\frac 1{r_{\min}}= \frac 1{r_1} -\frac {\rho_1}N$ and $\frac
1{r_{\max}} = \frac 1{r_1}
-\frac {\rho_1}{N-1}$.
If we take  $\rho_1 =-\left( N-\frac 32 \right)\left(\frac 12 - \frac
1{r_1}\right)+1$ in Theorem
\ref{str-in},   we have $\rho_1 \geq 0$ if and only if $  r_1 \leq
\frac{2(2N-3)}{2N-7}$.
Taking into account also the condition
 $r_1 \geq 2$, which corresponds to $\rho_1 \leq 1$, we obtain by
\eqref{sob}  the extremal spaces
\[
\begin{aligned}
u = v+w & \in  \spc{\infty}{I} 122 \cap \spc {2N-3} I 02{
\frac{2(2N-3)}{2N-7}}\\
& \subset \sp {\infty}{I} {L^{\frac  {2N}{N-2}}(\hn)\cap L^{\frac
{2(N-1)}{N-3}}(\hn)}
\cap \sp {2N-3} I{L^{ \frac{2(2N-3)}{2N-7}}(\hn)}.
\end{aligned}
\]
On the other hand, taking $\rho_1 =-\left( N-\frac 12 \right)\left(\frac 12
- \frac
1{r_1}\right)+1$ in Theorem \ref{str-in},   we have $\rho_1 \geq 0$  if and
only if $ r_1 \leq
\frac{2(2N-1)}{2N-5}$. The other bound
 $r_1 \geq 2$ still corresponds to $\rho_1 \leq 1$ and  we obtain therefore
the extremal spaces
\[
\begin{aligned}
u = v+w & \in  \spc
{\infty}{I} 122 \cap \spc {2N-1} I 02{ \frac{2(2N-1)}{2N-5}}\\
& \subset \sp {\infty}{I} {L^{\frac  {2N}{N-2}}(\hn)\cap L^{\frac
{2(N-1)}{N-3}}(\hn)}  \cap \sp
{2N-1} I{L^{ \frac{2(2N-1)}{2N-5}}(\hn)}.
\end{aligned}
\]

By interpolation we obtain
$u\in \sp {p}{I} {L^r(\hn)}$ with $0\leq \frac 2p \leq \frac 12 -\frac 1r$ and
$(N-1)\left(\frac 12 -\frac 1r\right) -1 \leq \frac 1p \leq N\left( \frac
12 -\frac 1r\right) -1$.

\hfill $\scriptstyle\bullet$\vskip\theorempostskipamount\par

\section{About the sharpness of the dispersive
estimates}\leavevmode\par\label{sharp}
We end up this paper by discussing the sharpness of the dispersive estimate
obtained in Proposition \ref{Dispe}.
Let us define the functions $v_j \in \srad$, $j \in \Z,$ by
\begin{equation*}
    \widehat{v_j}(m,\lambda)=       \left\{
         \begin{array}{ll}
             R(2^{-2j}(4n\lambda+\lambda^2)), & \hbox{if} \; m=0,\,\lambda>0,\\
             0, & \hbox{otherwise.} \\
         \end{array}
      \right.
\end{equation*}
\begin{lem}\label{Bes}
For any $\rho\in \R$ there exists $C_{\rho}>0$ such that
\begin{equation*}
    \|v_j\|_{\bL \rho11} \leq C_{\rho}2^{j\rho},\quad j\in\Z.
\end{equation*}
\end{lem}
\begin{proof}
We just have to prove the uniform estimate $\|v_j\|_{\lp1}\leq C,\;j \in
\Z$. Indeed
\begin{equation*}
    \widehat{v_j \ast \psi_k}(m,\lambda)=
           \left\{
             \begin{array}{ll}
               \RL j \RL k, & \hbox{if}\; m=0, \lambda >0\\
               0, & \hbox{otherwise} \\
             \end{array}
           \right.
\end{equation*}
implies $v_j\ast\psi_k=0$ if $|j-k|\geq2$. Therefore by \eqref{normauno}
\begin{equation*}
    \begin{split}
    \|v_j\|_{\bL \rho11} &
=\sum_{k=j-1}^{j+1}2^{k\rho}\|v_j\ast\psi_k\|_{\lp1}\\
        & \leq C_{\rho}\|v_j\|_{\lp1}2^{j\rho},
    \end{split}
\end{equation*}
where $C_{\rho}$ depends only on $\rho$.

Let us estimate $\|v_j\|_{\lp1}:$
\begin{equation}\label{a}
    |v_j(z,s)| =\cn \left|\i {e^{-i\lambda s}\RL j
e^{-\lambda|z|^2}\lambda^n}\right|,
\end{equation}
where $\lambda_1=\sqrt{4n^2+2^{2j-2}}-2n$ and
$\lambda_2=\sqrt{4n^2+2^{2j+2}}-2n$.
Then for $s \neq 0$
\begin{align}
    |v_j(z,s)|
        & =  \cn   \left|\i{(\frac{d}{d\lambda}e^{-i\lambda s})
                   \frac{\RL
je^{-\lambda|z|^2}\lambda^n}{is}}\right|\nonumber\\
        & = \cn  \frac{1}{|s|}\left|\i{e^{-i\lambda s}
                   \frac{d}{d\lambda}(\RL
je^{-\lambda|z|^2}\lambda^n)}\right|\nonumber\\
        &= \cn  \frac{1}{s^2}\left|\i{e^{-i\lambda s}
                   \frac{d^2}{d\lambda^2}(\RL
je^{-\lambda|z|^2}\lambda^n)}\right|.\label{b}
\end{align}
So we have two possible ways to estimate $|v_j(z,s)|$: using \eqref{a}
\begin{equation}\label{Stima1}
    |v_j(z,s)|\leq
         \left\{
          \begin{array}{ll}
             C2^{j(n+1)}e^{-\frac{2^j|z|^2}{C}}, & j\geq 0 \\
             C2^{2j(n+1)}e^{-\frac{2^{2j}|z|^2}{C}}, & j<0 \\
          \end{array}
         \right.
\end{equation}
or using \eqref{b}
\begin{equation}\label{Stima2}
    |v_j(z,s)|\leq
        \left\{
          \begin{array}{ll}
            \frac{C}{s^2}2^{j(n-1)}(1+2^j|z|^2+2^{2j}|z|^4)
            e^{-\frac{2^j|z|^2}{C}}, & j\geq0  \\
            \frac{C}{s^2}2^{2j(n-1)}(1+2^{2j}|z|^2+2^{4j}|z|^4)
            e^{-\frac{2^{2j}|z|^2}{C}}, & j<0 . \\
          \end{array}
        \right.
\end{equation}
For $j \geq0$, we have by \eqref{Stima1}
\begin{equation*}
    \begin{split}
      \int_{\{|s|<2^{-j},\,z\in \C^n\}}|v_j(z,s)|dzds
        & \leq C2^{j(n+1)}(\int_{\{|s|<2^{-j}\}}ds)(\int_{\C^n}
                e^{-\frac{2^j|z|^2}{C}}dz)\\
        & =
C2^{j}(\int_{\{|s|<2^{-j}\}}ds)(\int_{\C^n}e^{-\frac{|w|^2}{C}}dw)\leq C
   \end{split}
\end{equation*}
and by \eqref{Stima2}
\begin{equation*}
    \begin{split}
    \int_{\{|s|\geq2^{-j},\,z\in \C^n\}}|v_j(z,s)|dzds
        &\leq C2^{j(n-1)}(\int_{\{|s|\geq2^{-j}\}}\frac{1}{s^2}ds)(\int_{\C^n}
          (1+2^j|z|^2+2^{2j}|z|^4)e^{-\frac{2^j|z|^2}{C}}dz)\\
        & =C2^{-j}(\int_{\{|s|\geq2^{-j}\}}\frac{1}{s^2}ds)
           (\int_{\C^n}(1+|w|^2+|w|^4)e^{-\frac{|w|^2}{C}}dw)
           \leq C.
\end{split}
\end{equation*}
Therefore
\begin{equation*}
    \|v_j\|_{\lp1}=\int_{\{|s|<2^{-j},\,z\in \C^n\}}|v_j(z,s)|dzds +
    \int_{\{|s|\geq2^{-j},\,z\in \C^n\}}|v_j(z,s)|dzds \leq C.
\end{equation*}
Similarly, for $j <0$, we have
\begin{equation*}
    \|v_j\|_{\lp1}=\int_{\{|s|<2^{-2j},\,z\in \C^n\}}|v_j(z,s)|dzds +
    \int_{\{|s|\geq2^{-2j},\,z\in \C^n\}}|v_j(z,s)|dzds \leq C.
\end{equation*}
\end{proof}\par\medskip
By the definition of the functions $v_j$ we have
\begin{equation}
\cos
(t\sqrt\LL)v_j(0,\sigma_jt)=C2^{Nj}\int_0^{+\infty}e^{-it2^{2j}g_j(x)}h_j(x)
dx + C2^{Nj}\int_0^{+\infty}e^{-it2^{2j}\tilde g_j(x)}h_j(x)dx
\label{eitlvjsjt}\end{equation}
where $\sigma_j$ is a constant depending only on $j$,
$g_j=g_{j,\sigma_j,0}$ and $h_j=h_{j,0,0}$ are
the functions defined in
\eqref{gjsx} and \eqref{hzmx} respectively, and
\begin{equation}\label{gtilde}
    \tilde g_j(x)=\frac{1}{n}\left(\sigma_j x-\sqrt{2^{2-2j}n^2x+x^2}\right).
\end{equation}
\begin{lem}\label{inf}
For any $j \in \Z$ let $\eta_j$ be a function in $C^2 (\R)$  with supp
$\eta_j \subset [a_j,b_j]$
and let $\gamma_j$ be a real-valued function in $C^4 ([a_j,b_j])$ with
$\gamma'_j(x_j)=0$ for some
$x_j \in (a_j,b_j)$ and $\gamma_j''(x)\neq 0$ for any $x\in [a_j,b_j]$.
Therefore there exists $T_j>0$ such that
\begin{equation*}
    \left|\int_{a_j}^{b_j}e^{-it2^{2j}\gamma_j(x)}\eta_j(x)dx\right|\geq
    \frac{\sqrt \pi}{2}\t 2^{-j}|\gamma_j''(x_j)|^{-\mez}|\eta_j(x_j)|,
\quad t>T_j.
\end{equation*}
\end{lem}
\begin{proof}
It is not restrictive to suppose $\gamma_j(x_j)=0$ and $\eta_j(x_j)\neq 0$.
Let $\xi_j$ be the function defined by
\begin{equation*}
    \xi_j(x)=       \left\{
     \begin{array}{ll}
         -\sqrt{\frac{2\gamma_j(x)}{\gamma_j''(x_j)}}, & x\in [a_j,x_j] \\
         \sqrt{\frac{2\gamma_j(x)}{\gamma_j''(x_j)}}, & x\in (x_j,b_j]. \\
     \end{array}
       \right.
\end{equation*}
It is not hard to check that $\xi_j\in C^3([a_j,b_j]),\,\xi_j'>0$ on
$[a_j,b_j]$ and $\xi_j'(x_j)=1$.
Performing the change of variable $y=\xi_j(x)$
\begin{equation*}
    \int_{a_j}^{b_j}e^{-it2^{2j}\gamma_j(x)}\eta_j(x)dx=
\int_{\xi_j(a_j)}^{\xi_j(b_j)}e^{-it2^{2j}\frac{\gamma_j''(x_j)}{2}y^2}\Phi_j(y)
dy,
\end{equation*}
where $\Phi_j\in C^2$ and
$\Phi_j(y)=(\eta_j(\xi_j^{-1}(y))(\xi_j^{-1})'(y),$  supp
$\Phi_j\subset[\xi_j(a_j),\xi_j(b_j)]$ and $\Phi_j(0)=\eta_j(x_j)$.
We can write
\begin{equation*}
\int_{\xi_j(a_j)}^{\xi_j(b_j)}e^{-it2^{2j}\frac{\gamma_j''(x_j)}{2}y^2}\Phi_j(y)
dy=J_{j,t}+K_{j,t}
\end{equation*}
where
\begin{equation*}
    J_{j,t}=\int_{-\infty}^{+\infty} e^{-it2^{2j}
            \frac{\gamma_j''(x_j)}{2}y^2}e^{-y^2}\Phi_j(0)dy=
\frac{\sqrt \pi
\eta_j(x_j)}{\sqrt{|1+it2^{2j}\frac{\gamma_j''(x_j)}{2}}|}\,
            e^{-\frac{i}{2}\arctan(t2^{2j}\frac{\gamma_j''(x_j)}{2})}
\end{equation*}
and
\begin{equation*}
    \begin{split}
    K_{j,t} & =\int^{+\infty}_{-\infty}
e^{-it2^{2j}\frac{\gamma_j''(x_j)}{2}y^2}(\Phi_j(y)-e^{-y^2}\Phi_j(0))dy\\
            & =-\int^{+\infty}_{-\infty}
               \frac{d}{dy}(e^{-it2^{2j}\frac{\gamma_j''(x_j)}{2}y^2})
               \frac{\Phi_j(y)-e^{-y^2}\Phi_j(0)}{it2^{2j}\gamma_j''(x_j)y}dy\\
            &
=\frac{1}{it2^{2j}\gamma_j''(x_j)}\int^{+\infty}_{-\infty}e^{-it2^{2j}
\frac{\gamma_j''(x_j)}{2}y^2}\frac{d}{dy}(\frac{\Phi_j(y)-e^{-y^2}\Phi_j(0)}{y})
dy
.
\end{split}
\end{equation*}
Therefore
\begin{equation*}
    \left|\int_{a_j}^{b_j}e^{-it2^{2j}\gamma_j(x)}\eta_j(x)dx\right|
    \geq |J_{j,t}|\left|1-\frac{|K_{j,t}|}{|J_{j,t}|}\right|,
\end{equation*}
and, since $y\mapsto \frac{\Phi_j(y)-e^{-y^2}\Phi_j(0)}{y}$ is a function
in $C^1(\R)$ whose
derivative is in $L^1(\R)$, as is possible to verify by direct calculation,
we have
\begin{equation*}
\frac{|K_{j,t}|}{|J_{j,t}|}\leq\frac{\sqrt{|1+it\,2^{2j}\frac{\gamma_j''(x_j)}{2
}|}}
    {\sqrt \pi \, |\eta_j(x_j)|}\frac{C_j}{|t\,2^{2j}\gamma_j''(x_j)|}
    \leq C_j'\t \leq \mez, \quad t\geq 4(C_j')^2,
\end{equation*}
where $C_j$ and $C_j'$ are positive constants depending on $j$ but not on $t$.
Thus we obtain
\begin{equation*}
    \left|\int_{a_j}^{b_j}e^{-it2^{2j}\gamma_j(x)}\eta_j(x)dx\right|\geq
    \frac{\sqrt \pi}{2}\t2^{-j}|\gamma_j''(x_j)|^{-\mez}|\eta_j(x_j)|,
\quad t>T_j.
\end{equation*}
\end{proof}\par\medskip
Going back to (\ref{eitlvjsjt}), for any $j \in \Z$ we can fix $x_j>0$ such
that
$4x_j+\frac{2^{2j}x_j^2}{n^2}=1$ and $\sigma_j<0$ such that
$g_j'(x_j)=0$. By Lemma \ref{inf} and \eqref{baasta} we
obtain the following lower
estimates for $t>T_j$:
\begin{equation}\label{s3}
     \left|\int_0^{+\infty}e^{-it2^{2j}g_j(x)}h_j(x)dx\right|\geq
         \left\{
           \begin{array}{ll}
             C \t2^{-(n+\frac{1}{2})j}, & \hbox{if}\; j\geq 0 \\
             C \t2^{-\frac{j}{2}}, & \hbox{if}\; j< 0.\\
           \end{array}
         \right.
\end{equation}
In order to estimate the last integral in \eqref{eitlvjsjt} we first remark
that $\tilde g_j\,'(x)<0$ for any $x\in {\rm supp}\,h_j\subset[a_j,b_j]$.
Performing the change of variable $y=\tilde g_j(x)$
\begin{equation*}
    \left|\int_{a_j}^{b_j}e^{-it2^{2j}\tilde g_j(x)}h_j(x)dx\right|=
    \left|\int_{\tilde g_j(a_j)}^{\tilde g_j(b_j)}e^{-it2^{2j}y}H_j(y)dy\right|
\end{equation*}
where $H_j\in C^\infty$ and $H_j(y)=h_j(\tilde g_j\,^{-1}(y))(\tilde
g_j\,^{-1})(y)$, supp $H_j\subset [\tilde g_j(b_j),\tilde g_j(a_j)]$.
Then, for any $j\in \mathbb Z$ there exist $C_j, T'_j >0$ such that
\begin{equation}\label{s3a}
    \left|\int_0^{+\infty}e^{-it2^{2j}\tilde
g_j(x)}h_j(x)dx\right|=|\widehat H_j(t2^{2j})|
    \leq C_jt^{-1}, \quad t>T'_j.
\end{equation}
By \eqref{eitlvjsjt}, \eqref{s3} and \eqref{s3a} there exists $T''_j>0$
such that for $t>T''_j$:
\begin{equation}\label{s31}
     \|\cos (t\sqrt\LL)v_j\|_{L^\infty(\hn)}\geq
         \left\{
           \begin{array}{ll}
             C \t2^{(N-n-\frac{1}{2})j}, & \hbox{if}\; j\geq 0 \\
             C \t2^{(N-\frac{1}{2})j}, & \hbox{if}\; j< 0.\\
           \end{array}
         \right.
\end{equation}
\smallskip\\
\textit{Sharpness in $t$.}
Estimates (\ref{s31}) give for instance
\begin{equation*}
    \|\cos t\sqrt\LL v_0\|_{L^\infty(\hn)} \geq C\t, \quad t>T_0.
\end{equation*}
So the decay in $t$ in Proposition 1 cannot be improved.\smallskip\\
\textit{Sharpness in $\rho$.} Let us suppose that for some $\rho\in \R$ the
estimate
$\|\cos t\sqrt\LL f\|_{L^\infty(\hn)}
\leq C_{\rho} |t|^{-\mez} \|f\|_{\bL \rho11}$ holds for any $f\in\sh$. In
particular, by Lemma \ref{Bes},
$\|\cos t\sqrt\LL v_j\|_{L^\infty(\hn)} \leq C_{\rho} |t|^{-\mez}
2^{j\rho}, \; j
\in \Z$.
Estimates
\eqref{s31} force
$\rho\in[N-n-\frac{1}{2},
N-\frac{1}{2}]$.
\medskip\\
\textbf{Final remarks.} We would  like to emphasise that
there is no hope to
obtain a dispersive inequality as in Proposition 1
with the spaces
${\bD \rho qr}$.  Let us define the functions $w_j \in \srad$, $j \in \Z$, by
\begin{equation*}
    \widehat{w_j}(m,\lambda)=       \left\{
         \begin{array}{ll}
             R(2^{2-2j}n\lambda), & \hbox{if} \; m=0,\,\lambda>0\\
             0, & \hbox{otherwise.} \\
         \end{array}
      \right.
\end{equation*}
By the inversion formula \eqref{fx}
\begin{equation*}
    \begin{split}
    w_j(z,s) & =\cn\int_0^{+\infty} e^{-i\lambda
s}R(2^{2-2j}n\lambda)e^{-\lambda |z|^2}
                      \lambda^nd\lambda\\
        & =\cn 2^{Nj}\int_0^{+\infty} e^{-i\nu 2^{2j}s}R(4n\nu)e^{-\nu
|2^jz|^2}\nu^nd\nu\\
        & = 2^{Nj}w_0(2^jz,2^{2j}s).
\end{split}
\end{equation*}
Therefore $\|w_j\|_{\lp 1}=\|w_0\|_{\lp 1}$. This implies, as for the
functions $v_j$ (see the
proof of Lemma \ref{Bes}), that
$\|w_j\|_{\bD
\rho11}\leq C_\rho2^{j\rho}$, where $C_\rho$ depends only on $\rho$.
By the definition of $w_j$ we have
\begin{equation*}
\cos t\sqrt\LL
w_j(0,\sigma_jt)=C2^{Nj}\int_0^{+\infty}e^{-it2^{2j}g_j(x)}k(x)dx +
C2^{Nj}\int_0^{+\infty}e^{-it2^{2j}\tilde g_j(x)}k(x)dx
\end{equation*}
where $\sigma_j$ is a constant depending only on $j$,
$g_j=g_{j,\sigma_j,0}$ and $\tilde g_j$ are the functions defined in
(\ref{gjsx}) and \eqref{gtilde} respectively, and
$k(x)=R(4x)\frac{x^n}{n^{n+1}}$. For any $j\in\Z$ we fix
$x_j=\frac{1}{4}$ and $\sigma_j<0$ such that
$g_j'(\frac{1}{4})=0$. Arguing as before (see the proof of \eqref{s31}) we
obtain for $t>T_j$
\begin{equation*}
     \|\cos t\sqrt\LL w_j\|_{L^\infty(\hn)} \geq
         \left\{
           \begin{array}{ll}
             C \t2^{(N+1)j}, & \hbox{if}\; j\geq 0\\
             C \t2^{(N-\frac{1}{2})j}, & \hbox{if}\; j< 0.\\
           \end{array}
         \right.
\end{equation*}
These estimates imply that there is no $\rho\in\R$ for which
$\|\cos t\sqrt\LL f\|_{L^\infty(\hn)}
\leq C |t|^{-\mez}\|f\|_{\bD \rho11}$ for any $f\in\sh$.
\medskip

As a conclusion we would like to remark that analysing the wave equation
related to the Kohn-Laplacian $\Delta$ with the spaces $\bpl{\rho}{q}{r}$
we obtain the dispersive inequality for the wave semigroup: for any
$\rho\in[N-\frac{3}{2},N-\mez]$
\begin{equation*}
\left\|e^{-it\sqrt{\Delta}}f \right\|_{\lp \infty}\leq C_{\rho}
|t|^{-\mez}\|f\|_{\bpl{\rho}{1}{1}},\qquad f\in\sh,\ t\in\R^{*}.
\end{equation*}
This result does not give Proposition 1 (unless $u_1=0$) because estimate
\eqref{fmv} does not hold with $L=\Delta$ and $\Delta_ju=u\ast\psi_j$.

Finally, for the Schr\"odinger equation
related to the full
Laplacian, by Proposition \ref{inclusioni} and \cite[Corollary 10]{FV} we
have the dispersive
estimate
\begin{align}
\left\|e^{-it\LL}f \right\|&_{\lp \infty}\leq C
t^{-\mez}\|f\|_{\bpl{N-2}{1}{1}},\qquad f\in\sh,\ t>0.\label{dispSch}
\end{align}
By a direct computation as in Section 4 the estimate \eqref{dispSch} cannot be
improved. So
the behaviour of the Schr\"odinger operator $e^{-it\LL}$ by analysing it
with the spaces
$\bpl{\rho}{q}{r}$ is the same as in [FV] with the spaces $\bpd{\rho}{q}{r}$.

\end{document}